\documentclass{amsart}
\usepackage{amsmath}
\usepackage{url}

\theoremstyle{plain}
\newtheorem{theorem}{Theorem}

\newtheorem{lem}[theorem]{Lemma}

\newtheorem{cor}[theorem]{Corollary}

\newtheorem *{Lambek-Moser Theorem}{Lambek-Moser Theorem}

\theoremstyle{remark}

\theoremstyle{definition}
\newtheorem{definition}[theorem]{Definition}
\begin{document}

\title{On the Lambek-Moser Theorem}
\author{Yuval Ginosar }
\address{Department of Mathematics, University of Haifa, Haifa 31905, Israel}
\email{ginosar@math.haifa.ac.il}
\date{\today}
\maketitle

\begin{abstract}
We suggest an alternative proof of a theorem due to Lambek and
Moser using a perceptible model.
\end{abstract}
\section{}
The
notion of invertibility of sequences which take their values in
$\mathbb{Z}^+\cup \{\infty\}$ (where $\mathbb{Z}^+$ denotes the set of non-negative integers) was
introduced by J. Lambek and L. Moser. Adopting their terminology \cite{LM},
such sequences are called sequences of {\it numbers}.
\begin{definition}
Two sequences
$\bar{f}=(f(n))_{n=1}^{\infty}, \bar{g}=(g(n))_{n=1}^{\infty}$ of numbers are
{\it mutually inverse} if for every $m,n\in \mathbb{Z}^+$ either
$f(m)<n$ or $g(n)<m$, but not both.
\end{definition}
It is shown \cite[Theorem 1]{LM} that a sequence of numbers $(f(n))_{n=1}^{\infty}$
has an inverse if and only if it is non-decreasing.
In this case, the unique inverse $(g(n))_{n=1}^{\infty}$ is given by
\begin{equation}\label{inv}
g(n)=| \{m| f(m)<n\}|.
\end{equation}
It follows that the inverse of $(g(n))_{n=1}^{\infty}$ is again
 $(f(n))_{n=1}^{\infty}$.

Any non-decreasing sequence of numbers $\bar{f}=(f(n))_{n=1}^{\infty}$, determines a set of positive integers
$$\hat{f}:=\{ n+f(n) \}_{f(n)\in\mathbb{Z}^+}.$$
The correspondence $\bar{f}\mapsto \hat{f}$ between the non-decreasing sequences of numbers and
the sets of positive integers is one-to-one.
The following partitioning theorem is established.
\begin{Lambek-Moser Theorem} \cite[Theorem 2]{LM}
Two non-decreasing sequences of numbers $\bar{f}=(f(n))_{n=1}^{\infty}, \bar{g}=(g(n))_{n=1}^{\infty}$
are mutually inverse if and only if the sets $\hat{f}$ and $\hat{g}$
are complementary, that is they disjointly cover the positive integers.
\end{Lambek-Moser Theorem}
In this note we suggest an alternative proof of the Lambek-Moser theorem, by applying the
running model which was introduced in \cite{GY}.
Another visual proof was given by E.W. Dijkstra \cite{D}.
The reader is referred to \cite{ZWS} for a detailed bibliography on complementary sequences and related topics.

\section{}
\noindent {\bf The Model.} Let $X$ and $Y$ be two athletes running around a circular stadium in opposite directions,
starting at time $t=0$. 
Each time one of these athletes crosses a fixed point $\mathcal{O}$, the number of
their meetings since the time $t=0$ (if a meeting has already happened) is recorded for this athlete.
Notice that the distance traveled by both athletes between two consecutive meetings is precisely the length
of the stadium.
Now, since the athletes meet countably many times, we choose the point $\mathcal{O}$ in which none of these meetings occurs.
Then it is clear that between two consecutive meetings,
exactly one of the two of them crosses $\mathcal{O}$.
As a result, the set $\mathcal{S}_X$ recorded for $X$ and the
set $\mathcal{S}_Y$ recorded for $Y$ partition the set of positive integers.
Conversely, it is also necessary for these sets to partition the set of positive integers that
none of their meetings occurs at $\mathcal{O}$.

\section{}
Normalize the length of the stadium to be 1, place $\mathcal{O}$ in 0(mod 1), and let
$$\begin{array}{rcl}
\varphi:\mathbb{R}^+&\rightarrow &\mathbb{R}^+\\
t&\mapsto &\varphi(t)
\end{array}$$
be a strictly increasing continuous time function describing the motion of $X$.
Let $\psi(t)=t$ be the motion function of $Y$, who is running in the opposite direction.
Then $Y$ crosses $\mathcal{O}$ exactly in integer time units. Since the relative motion function of $X$ and $Y$
is $\varphi(t)+t$, and since together they travel a unit between two consecutive meetings,
the number of times $X$ and $Y$ meet until time $t$
is $\lfloor \varphi(t)+t\rfloor$, where $\lfloor \cdot\rfloor$ is the floor integer part function.
Therefore, the set $\mathcal{S}_Y$ of positive integers recorded for $Y$ satisfies
$$\mathcal{S}_Y=\{\lfloor \varphi(n)+n\rfloor\}_{n\in\mathbb{Z}^+}.$$

Next, $X$ crosses the point $\mathcal{O}$ each and every time $\varphi(t)\in\mathbb{Z}^+$.
Thus, the set $\mathcal{S}_X$ recorded for $X$ satisfies
$$\mathcal{S}_X=\{\lfloor \varphi(t)+t\rfloor\}_{\varphi(t)\in\mathbb{Z}^+}.$$
We can describe $\mathcal{S}_X$ in another way.
Since $\varphi$ is continuous and strictly increasing, it maps $\mathbb{R}^+=(0,\infty)$ onto an open segment
$I:=(0,M)$ (where $0<M\leq \infty$), and admits an increasing, continuous inverse
$\varphi^{-1}:I\to\mathbb{R}^+$.
Then $$\mathcal{S}_X=\{\lfloor n+\varphi^{-1}(n)\rfloor\}_{n\in\mathbb{Z}^+\cap I}.$$
By the above argument, the sets $\mathcal{S}_X$ and $\mathcal{S}_Y$
partition the positive integers if and only if $X$ and $Y$ never meet at $\mathcal{O}$.
But $X$ and $Y$ do meet at $\mathcal{O}$ at time $t$ exactly if both $t$ and $\varphi(t)$ are in $\mathbb{Z}^+$.
We obtain
\begin{cor}\label{cont}
Let $\varphi:\mathbb{R}^+\rightarrow \mathbb{R}^+$ be a strictly increasing continuous function and let
$\varphi^{-1}:\rm {Im}(\varphi)\rightarrow \mathbb{R}^+$ be its inverse. Then the sets
$\{\lfloor \varphi(n)+n\rfloor\}_{n\in\mathbb{Z}^+}$
and $\{\lfloor n+\varphi^{-1}(n)\rfloor\}_{n\in\mathbb{Z}^+\cap\rm{Im}(\varphi)}$
partition the positive integers if and only if $\varphi(\mathbb{Z}^+)\cap \mathbb{Z}^+= \emptyset$.
\end{cor}

\section{}\label{prf}
In order to exploit Corollary \ref{cont} to prove the Lambek-Moser Theorem, we need two lemmas.
The first observation is easily verified by distinguishing between three types of sequences (see \cite[\S 2]{LM}).
\begin{lem}\label{int}
Let $(f(n))_{n=1}^{\infty}$ and $(g(n))_{n=1}^{\infty}$ be mutually inverse sequences of numbers.
Then at least one of these sequences does not admit $\infty$ as a value, in other words, it is a sequence of (non-negative) integers.
\end{lem}
The second lemma is straightforward:
\begin{lem}\label{cond}
Let  $(f(n))_{n=1}^{\infty}$ be a non-decreasing sequence of (non-negative) integers.
Then there exists a strictly increasing continuous function
$\varphi:\mathbb{R}^+\rightarrow \mathbb{R}^+$ such that
$\lfloor \varphi(n)\rfloor=f(n)$, for every $n\in \mathbb{Z}^+$.
Moreover, $\varphi$ can be chosen such that
\begin{equation}\label{cap}
\varphi(\mathbb{Z}^+)\cap \mathbb{Z}^+=\emptyset.
\end{equation}
\end{lem}
\noindent {\bf Proof of the Lambek-Moser Theorem.}
By Lemma \ref{int} we may assume that $(f(n))_{n=1}^{\infty}$
is sequence of non-negative integers (else $(g(n))_{n=1}^{\infty}$ is).
Next, by Lemma \ref{cond}, there exists
a strictly increasing continuous function $\varphi:\mathbb{R}^+\rightarrow \mathbb{R}^+$ such that
for every $n\in\mathbb{Z}^+$, both\\
(a) $\lfloor\varphi(n)\rfloor=f(n)$, and\\
(b) $\varphi(n)\notin \mathbb{Z}^+$.

We can now use Corollary \ref{cont} to deduce that the sets
$$\{f(n)+n\}_{n=1}^{\infty}=\{\lfloor \varphi(n)+n\rfloor\}_{n\in\mathbb{Z}^+}$$ and
$$\{\lfloor \varphi^{-1}(n)+n\rfloor\}_{n\in\mathbb{Z}^+\cap\text {Im}(\varphi)}$$ are complementary.
By the conditions on $\varphi$, using the alternative characterization \eqref{inv} for the inverse sequence,
the inverse of $\bar{f}=(f(n))_{n=1}^{\infty}=(\lfloor \varphi(n)\rfloor)_{n=1}^{\infty}$
is $\bar{g}$ given by
\begin{equation}\label{chareg}
g(n)=\left \{
\begin{array}{cl}
\lfloor \varphi^{-1}(n)\rfloor & \text { if $n\in$ Im($\varphi$)}\\
\infty & \text { otherwise.}
\end{array}
\right.
\end{equation}
Consequently, the complement of $\hat{f}$ in $\mathbb{Z}^+$ is $\hat{g}$.
Since the correspondence $\bar{f}\mapsto \hat{f}$ is one-to-one,
the proof of the theorem is complete.  $\square$\\
{\bf Remark.}
Note that S. Beatty's celebrated theorem \cite{BSOH} follows from Corollary \ref{cont} by taking
$\varphi(t):=\lambda\cdot t$,
where $\lambda>0$ (and then $\varphi^{-1}(t)=\frac{1}{\lambda}\cdot t$), that is the case where the speeds of both
athletes are constant.

\noindent{\bf Acknowledgement.}
Many thanks to my colleagues G. Moran and D. Blanc for their useful suggestions.

\end{document}